\documentclass{amsart}%
\usepackage{amsfonts}
\usepackage{amsmath}
\usepackage{amssymb}
\usepackage{graphicx}%
\providecommand{\U}[1]{\protect\rule{.1in}{.1in}}
\newtheorem{theorem}{Theorem}
\theoremstyle{plain}
\newtheorem{acknowledgement}{Acknowledgement}

\newtheorem{corollary}{Corollary}

\newtheorem{definition}{Definition}

\newtheorem{lemma}{Lemma}

\newtheorem{proposition}{Proposition}
\newtheorem{remark}{Remark}

\numberwithin{equation}{section}
\begin{document}
\title[CY\ Manifolds with Locally Symmetric Moduli Space]{Calabi-Yau Manifolds whose Moduli Spaces are Locally Symmetric Manifolds and
No Quantum Corrections}
\author{Andrey N. Todorov}
\address{UCSC, Department of Mathematics, Santa Cruz, CA 95064\\
Bulgarian Academy of Sciences, Institute of Mathematics, Sofia 1113 Bulgaria}
\email{todorov.andrey@gmail.com}
\date{June 26, 2008}

\begin{abstract}
It is observed that there are a natural sequence of CY manifolds M$_{g}$ that
are double covers of $\mathbb{CP}^{g}$ ramified over $2g+2$ hyperplanes and
some of them are obtained from the Jacobian $J(C_{g})$ of hyper-elliptic curve
$C_{g}$ by an action of an involution of a semi-direct product of $g-1$ copies
of groups of order two with the symmetric group of g elements. This
construction generalizes Kummer surfaces. We observed that the $g^{2}$
Kodaira-Spencer classes on the Jacobian are invariant under the action of the
above group. They form a basis of $H^{1}\left(  \text{M}_{g},T_{\text{M}_{n}%
}^{1,0}\right)  .$ Since on complex torus we have $\left[  \phi_{1},\phi
_{2}\right]  =0,$ then for any Kodaira-Spencer classes $\phi_{1}$ and
$\phi_{2}$ on these special CY manifolds, we have $\left[  \phi_{1},\phi
_{2}\right]  =0.$ From here it follows that the moduli space of those CY
manifolds is a locally symmetric space $\Gamma\backslash\mathbb{SU}%
(n,n)/\mathbb{S}(\mathbb{U}(n)\times\mathbb{U}(n)).$

\end{abstract}
\maketitle
\tableofcontents

\section{Introduction}

\subsection{General Remarks}

In the last fifteen years CY manifolds attracted the attention of
mathematicians and physicists. The influence of the ideas from String Theory
on certain fields of mathematics were amazing. Some of the ideas led to the
solutions of some classical problems. One such example is the applications of
the ideas of mirror symmetry to enumerative geometry and can be considered as
solution of one of the Hilbert's Problems.

The motivation of writing this paper is the recent paper \cite{GV} of S. Gukov
and C. Vafa. In \cite{GV} a problem to find CY manifolds whose moduli space
contains everywhere dense subset of CY manifolds of CM type was posed. This
problem is connected with the Conformal Field Theories. The Rational Conformal
Field Theory was introduced in \cite{FQS}. See also \cite{MS} and \cite{MS1}.
It is expected that sigma models with targets CY manifold of CM type will give
a Rational Conformal Field Theory. It was pointed out that the existence of
everywhere dense subset of RCFT may lead to a better understanding of the
Conformal Field Theory. According to a conjecture of F. Oort and Y. Andr\'{e},
the points in the moduli space of polarized algebraic manifold corresponding
to manifolds of CM types is dense only when the moduli space is a locally
symmetric space. See \cite{A1}, and \cite{O}. For some interesting relations
between the arithmetic of CY manifolds and Rational Conformal Field Theory see
\cite{Sch}. According to Satake $\Gamma\backslash\mathfrak{h}_{g,g}$ is a
Shimura variety. See \cite{Sat}. It is natural to expect that the points of
$\Gamma\backslash\mathbb{SU}(g,g)/\mathbb{S}(\mathbb{U}(g)\times
\mathbb{U}(g))$ that corresponds to abelian varieties of CM type will
correspond to CY manifolds of CM type. This question will be addressed in
another paper.

The moduli spaces of Riemann surfaces of genus at least two are not locally
symmetric spaces. It seems that it is very rare that the moduli spaces of
projective algebraic varieties is a locally symmetric space. For example the
moduli space of one and two dimensional projective varieties with zero
canonical class are locally symmetric spaces. In dimension one the analogue of
CY manifold is the elliptic curve, in dimension two it is the K3 surface. The
moduli space of elliptic curves is $\Gamma_{E}\backslash\mathbb{SU}%
(1,1)/\mathbb{S}(\mathbb{U}(1)\times\mathbb{U}(1)).$ The moduli space of
algebraic polarized K3 surfaces is the locally symmetric space $\Gamma
_{K3,L}\backslash\mathbb{SO}_{0}(2,19)/\mathbb{SO}(2)\times\mathbb{SO}(19),$
where $\Gamma_{E}$ and $\Gamma_{K3,L}$ are arithmetic groups. The moduli space
of abelian varieties are also locally symmetric spaces.

Viehweg proved that the moduli space of polarized CY manifolds exists and it
is a quasi-projective variety. See \cite{vw}. In general the moduli space of
polarized CY manifolds is not a locally symmetric space. It seems that up to
now only one example of a CY manifold is known whose moduli space is a locally
symmetric space. The example is the Borcea-Voisin manifold which is obtained
from a K3 surface $X$ with an involution $\sigma_{1}$ without fixed points on
$X$ and elliptic curve $E$ as follows: $X\times E/\left(  \sigma_{1}%
\times\sigma_{2}\right)  ,$ where $\sigma_{2}$ is the involution on $E.$ Its
moduli space is a product of two locally symmetric spaces.

Dick Gross classified all possible symmetric spaces that are tube domains
which can be period domains of Variations of Hodge Structures of CY
threefolds. See \cite{Gr}. He posed a very interesting question to find CY
threefolds, whose moduli spaces are the list presented in \cite{Gr}.

There is a natural sequence of CY manifolds that are obtained as double covers
of $\mathbb{CP}^{g}$ ramified over $2g+2$ hyperplanes in general position.
Elliptic curves are double covers of $\mathbb{CP}^{1}$ ramified over four
points in general position. Double covers of $\mathbb{CP}^{2}$ ramified over
six lines in general position are K3 surfaces with $15$ double points whose
moduli space is isomorphic to $\Gamma\backslash\mathbb{SU}(2,2)/\mathbb{S}%
(\mathbb{U}(2)\times\mathbb{U}(2)).$ In this paper we prove that moduli space
of non singular CY manifolds M$_{g}$ that are obtained from the resolution of
the singularities of the double covers of $\mathbb{CP}^{g}$ ramified over
$2g+2$ hyperplanes in general position is the locally symmetric space
$\Gamma\left\backslash \mathbb{SU}(g,g)/\mathbb{S}(\mathbb{U}(g)\times
\mathbb{U}(g))\right.  =\Gamma\backslash\mathfrak{h}_{g,g},$ where $\Gamma$ is
a subgroup in the mapping class group which preserve the polarization class.

Dolgachev conjectured that the moduli space of CY manifolds that are double
covers of $\mathbb{CP}^{g}$ ramified over $2g+2$ planes is the tube domain
$S_{g}(\mathbb{C})+iS_{g}^{+}(\mathbb{C}),$ where $S_{g}(\mathbb{C})$ is the
space of $g\times g$ Hermitian matrices and $S_{g}^{+}(\mathbb{C})$ is the
space of positive Hermitian matrices. So $S_{g}(\mathbb{C})+iS_{g}%
^{+}(\mathbb{C})$ is the tube domain realization of the symmetric space
$\mathfrak{h}_{g,g}=\mathbb{SU}(g,g)/\mathbb{S}(\mathbb{U}(g)\times
\mathbb{U}(g)).$

\subsection{The Description of the Ideas of the Proofs of the Main Theorems}

The main observation of this paper is the following:

\textbf{Theorem 11. }\textit{Let M}$_{g}$\textit{ be a CY manifold obtained
from the resolution of singularities of a double cover of }$\mathbb{CP}^{g}%
$\textit{ ramified over }$2g+2$\textit{ hyperplanes in general position. Let
}$\phi_{1}$\textit{ and }$\phi_{2}\in H^{1}\left(  \text{M}_{g},T_{\text{M}%
_{g}}^{1,0}\right)  $\textit{ be two Kodaira-Spencer classes. Then }$\left[
\phi_{1},\phi_{2}\right]  =0.$

Theorem \textbf{11} follows from the following observation: Let $C_{g}$ be a
hyper-elliptic curve of genus $g>2$ such that the projection $\pi_{g}%
:C_{g}\rightarrow\mathbb{CP}^{1}$ has the ramification points $\lambda
_{1},...,\lambda_{2g+1}$ and $\infty.$ We can identify $J(C_{g})$ with
$Pic^{g-1}(C).$ The $\theta-$divisor of the Jacobian of a hyper-elliptic curve
can be identified with the image of $S^{g-1}(C_{g}).$ Notice that $S^{g}%
(C_{g})$ is obtained from $J(C_{g})$ by blowing up the subvariety of
codimension two which is the image of $S^{g-2}(C_{g})+2\lambda_{i}.$ Let us
consider $2g+2$ $\theta_{i}-$divisors obtained from $\theta$ by translates by
the $2n+1$points $\lambda_{i}-\infty.$ Notice that the $2g+2$ $\theta
-$divisors are tangent on $S^{g}(C_{g})$ to the blown up $S^{g-2}%
(C_{g})+2\lambda_{i}/$

We observed that on $\underset{g}{\underbrace{C_{g}\times...\times C_{g}}}$
the subgroup $\left(  \left.  \mathbb{Z}\right/  2\mathbb{Z}\right)
^{g-1}\ltimes S_{g}$ embedded into $\left(  \left.  \mathbb{Z}\right/
2\mathbb{Z}\right)  ^{g}\ltimes S_{g}$ as the sum of the order two elements to
be zero acts. It is easy to see that the quotient is a double cover of
$\mathbb{CP}^{g}$ with ramification divisors the $2g+2$ hyperplanes which are
the images of the $\theta-$divisors and are tangent to the image of the blown
up $S^{g-2}(C_{g})+2\lambda_{i}.$We observed that each vector of the basis
$\left\{  \overline{dz^{j}}\otimes\frac{\partial}{\partial z^{i}}\right\}  $
of $H^{1}\left(  J(C_{g}),T_{J(C_{g})}^{1,0}\right)  $ is invariant under the
action of $\left(  \left.  \mathbb{Z}\right/  2\mathbb{Z}\right)
^{g-1}\ltimes S_{g}.$ Since the fixed point set of $\tau_{g}$ has a
codimension two, then by Hartog's Theorem we can extend $\left\{  \left(
\nu_{n}\right)  _{\ast}\left(  \overline{dz^{j}}\otimes\frac{\partial
}{\partial z^{i}}\right)  \right\}  $ to a basis of $H^{1}\left(  \text{M}%
_{g},T_{\text{M}_{g}}^{1,0}\right)  .$ Thus each Kodaira Spencer class $\phi$
on M$_{g}$ can be written as follows:%
\[
\phi=%
{\displaystyle\sum\limits_{i,j=1}^{g}}
a_{\overline{j}}^{i}\overline{dz^{j}}\otimes\frac{\partial}{\partial z^{i}},
\]
where $a_{\overline{j}}^{i}$ are constant. The definition of a $\left[
\phi_{1},\phi_{2}\right]  $ implies that $\left[  \phi_{1},\phi_{2}\right]
=0$ for any $\phi_{k}\in H^{1}\left(  \text{M}_{g},T_{\text{M}_{g}}%
^{1,0}\right)  .$ This fact implies that the moduli space of the CY manifolds
that are double covers of $\mathbb{CP}^{g}$ are locally symmetric spaces. We
will show in the next paper and there are no quantum corrections to the Yukawa couplings.

\begin{acknowledgement}
I want to express my deep gratitude toward Igor Dolgachev for his remarks. I
disagree with some of his suggestions and comments about other papers related
to this one. \textit{Special thanks to Maxim Kontsevich for his useful
remarks. The author is grateful to Cumrum Vafa for suggesting to look for an
abelian variety that covers the double cover of }$\mathbb{CP}^{g}$
\textit{ramified over }$2g+2$ \textit{hyperplanes in general position. His
suggestion played an important role in this paper. The author wants to thank
Professor Yau, Professor Kefeng Liu and Professor Lu for their useful and
helpful comments and numerous discussions on this topic. Special thanks to
Sasha Goncharov and Yu. Manin for their interest and encouragement. }
\end{acknowledgement}

\section{Moduli Spaces of CY Manifolds}

\subsection{Local Deformation Theory of CY Manifolds}

\begin{definition}
\label{BELT}Let M be a complex manifold. Let
\[
\phi\in C^{\infty}\left(  \text{M},Hom\left(  \Omega_{\text{M}}^{1,0}%
,\Omega_{\text{M}}^{0,1}\right)  \right)  =C^{\infty}\left(  \text{M}%
,T_{\text{M}}^{1,0}\otimes\Omega_{\text{M}}^{0,1}\right)
\]
then we will call $\phi$ a Beltrami differential.
\end{definition}

The Beltrami differential \ $\phi$\ defines an integrable complex structure on
M if and only if the following equation holds: $\overline{\partial}\phi
=\frac{1}{2}\left[  \phi,\phi\right]  ,$ where%
\[
\left[  \phi,\phi\right]  |_{D^{n}}:=
\]%
\[
\sum_{\nu=1}^{n}\sum_{1\leqq\alpha<\beta\leqq n}\left(  \sum_{\mu=1}%
^{n}\left(  \phi_{\overline{\alpha}}^{\mu}\left(  \partial_{\mu}%
\phi_{\overline{\beta}}^{\nu}\right)  -\phi_{\overline{\beta}}^{\mu}\left(
\partial_{\nu}\phi_{\overline{\alpha}}^{\nu}\right)  \right)  \right)
\overline{dz}^{\alpha}\wedge\overline{dz}^{\beta}\otimes\frac{\partial
}{dz^{\nu}}.
\]

The main results in \cite{to89} are the following Theorems:

\begin{theorem}
\label{tod1}Let M be a CY manifold and let $\left\{  \phi_{i}\right\}  $ be a
basis of harmonic $(0,1)$ forms with coefficients in $T^{1,0}$ of
$\mathbb{H}^{1}($M$,T^{1,0}),$ then the equation $\overline{\partial}%
\phi=\frac{1}{2}\left[  \phi,\phi\right]  $ has a solution in the form:
$\phi(\tau)=\sum_{i=1}^{N}\phi_{i}\tau^{i}+\sum_{|I_{N}|\geqq2}\phi_{I_{N}%
}\tau^{I_{N}}$ and $\partial^{\ast}\left(  \phi(\tau)\right)  =0,$
\textit{where } $I_{N}=(i_{1},..,i_{N})$\ \ \textit{is a multi-index},
\ $\phi_{I_{N}}\in C^{\infty}($M$,\Omega^{0,1}\otimes T^{1,0}),$ $\phi_{I_{N}%
}\lrcorner\omega_{0}=\partial\psi_{I_{N}}$ $\tau^{I_{N}}=(\tau^{1})^{i_{1}%
}..(\tau^{N})^{i_{N}}$\textit{and there exists} $\varepsilon>0$ \textit{such
that } $\phi(\tau)\in C^{\infty}($M$,\Omega^{0,1}\otimes T^{1,0})$ for
$|\tau^{i}|<\varepsilon$ \ for $i=1,..,N.$
\end{theorem}

In \cite{to89} we proved the following basic fact:

\begin{theorem}
\label{tod2}Let $\omega_{0}$ \ be a holomorphic n-form on the n dimensional CY
manifold M. Let \ $\left\{  U_{i}\right\}  $ be a covering of M and let
$\left\{  z_{1}^{i},..,z_{n}^{i}\right\}  $ be local coordinates in $U_{i}$
such that $\omega_{0}|_{U_{i}}=dz_{1}^{i}\wedge...\wedge dz_{n}^{i}.$ Then for
each $\tau=(\tau^{1},..,\tau^{N})$ such that $|\tau_{i}|<\varepsilon$ the
forms on M defined as:
\[
\left.  \omega_{\tau}\right\vert _{U_{i}}:=(dz_{1}^{i}+\phi(\tau)(dz_{1}%
^{i}))\wedge..\wedge(dz_{n}^{i}+\phi(\tau)(dz_{n}^{i}))
\]
are globally defined complex n forms $\omega_{\tau}$ on M and, moreover,
$\ \omega_{\tau}$\ are closed holomorphic n forms with respect to the complex
structure on M defined by $\phi(\tau).$
\end{theorem}

\begin{corollary}
\label{tod3}We have the following Taylor expansion for the form $\omega_{\tau
}:$%
\begin{equation}
\left.  \omega_{\tau}\right\vert _{U_{i}}=\omega_{0}+\sum_{k=1}^{n}%
(-1)^{\frac{k(k-1)}{2}}\left(  \left(  \wedge^{k}\phi(\tau)\right)
\lrcorner\omega_{0}\right)  . \label{forms}%
\end{equation}
$\phi(\tau)$ is defined by Theorem \ref{tod1}$.$ (See \cite{to89}.)
\end{corollary}

\begin{corollary}
\label{tod4}We have%
\[
\left[  \omega_{\tau}\right]  =[\omega_{0}]+\sum[\left(  \omega_{0}%
\lrcorner\phi_{i}\right)  ]\tau^{i}-\sum[(\omega_{0}\lrcorner\left(  \phi
_{i}\wedge\phi_{j}\right)  )]\tau^{i}\tau^{j}+
\]
\begin{equation}
+\sum_{k\geq3}^{n}%
{\displaystyle\sum\limits_{0\leq i_{1}\leq...i_{k}\leq n_{i}}}
[\omega_{0}\lrcorner\left(  \phi_{i_{1}}(\tau)\right)  \wedge...\phi_{i_{k}%
}(\tau)], \label{tod41}%
\end{equation}
where $[\omega_{0}(n-k,k)]$ are forms of type $(n-k,k).$ (See \cite{to89}.)
\end{corollary}

\begin{definition}
\label{coor}We call this coordinate system a flat coordinate system. Let
$\mathcal{K\subset}\mathbb{C}^{N}$ be the polydisk defined by $|\tau
^{i}|<\varepsilon$ for every $1\leq i\leq N$, where $\varepsilon$ is chosen
such that for every $\tau\in\mathcal{K}$ , $\phi(\tau)\in C^{\infty}%
(M,\Omega^{0,1}\otimes T^{1,0}),$ where $\phi(\tau)$ is defined in \ Theorem
\ref{tod1}. On the trivial $C^{\infty}$ family M$\times\mathcal{K},$ we will
define for each $\tau\in\mathcal{K}$ an integrable complex structure
I$_{\phi(\tau)}$ on the fibre over $\tau$ of the family M$\times\mathcal{K}.$
So a complex analytic family $\pi:\mathcal{X\rightarrow K}$ of CY manifold was
defined. We will call this family the Kuranishi family. We introduce a
coordinate system $(\tau^{1},..,\tau^{N})$ in $\mathcal{K}$ as in Theorem
\ref{tod1}.
\end{definition}

\subsection{Global Theory of the Moduli Space of Polarized CY Manifolds}

\begin{definition}
We will call a pair (M; $\gamma_{1},...,\gamma_{b_{n}}$) a marked CY manifold
where M is a CY manifold and $\{\gamma_{1},...,\gamma_{b_{n}}\}$ is a basis of
$H_{n}$(M,$\mathbb{Z}$)/Tor.
\end{definition}

We proved in \cite{LTYZ} the existence of the Teichm\"{u}ller space:

\begin{theorem}
\label{teich}There exists a family $\mathcal{Z}_{L}\mathcal{\rightarrow
}\mathfrak{T}($M$),$ of marked polarized CY manifolds which possesses the
following properties: \textbf{a)} It is effectively parametrized, \textbf{b)
}For any marked CY manifold M of fixed topological type for which the
polarization class $L$ defines an imbedding into a projective space
$\mathbb{CP}^{N},$ there exists an isomorphism of it (as a marked CY manifold)
with a fibre M$_{s}$ of the family $\mathcal{Z}_{L}.$ \textbf{c) }The base has
dimension $h^{n-1,1}.$
\end{theorem}

\begin{corollary}
\label{teich1}Let $\mathcal{Y\rightarrow}\mathfrak{K}$ be any family of marked
CY manifolds, then there exists a unique holomorphic map $\phi:\mathfrak{K}%
\rightarrow\mathfrak{T}($M$)$ up to a biholomorphic map $\psi$ of M which
induces the identity map on $H_{n}($M$,\mathbb{Z}).$
\end{corollary}

\begin{corollary}
\label{teich2}The Teichm\"{u}ller space $\mathfrak{T}($M$)$ has finite number
of irreducible components. From now on we will denote by $\mathcal{T}$(M) the
irreducible component of the Teichm\"{u}ller space $\mathfrak{T}($M$)$ that
contains our fixed CY manifold M.
\end{corollary}

\begin{definition}
\label{mapgr}We will define the mapping class group $\Gamma$(M) of any compact
C$^{\infty}$ manifold M as follows: $\Gamma($M$)=Diff_{+}($M$)/Diff_{0}($M$),$
where $Diff_{+}($M$)$ is the group of diffeomorphisms of M preserving the
orientation of M and $Diff_{0}($M$)$ is the group of diffeomorphisms isotopic
to identity. Let $L\in H^{2}($M$,\mathbb{Z})$ be the imaginary part of a
K\"{a}hler metric. Let $\Gamma_{L}:=\{\phi\in\Gamma$(M)$|\phi(L)=L\}.$
\end{definition}

It is a well know fact that the moduli space of polarized algebraic manifolds
$\mathcal{M}_{\omega}($M$)=\mathfrak{T}($M$)/\Gamma_{\omega}.$ In \cite{LTYZ}
the following Theorem was proved.

\begin{theorem}
\label{Vie}There exists a subgroup of finite index $\Gamma$ of $\ \Gamma_{L}$
such that $\Gamma$ acts freely on the component T(M) of the Teichm\"{u}ller
space $\mathfrak{T}($M$),$ stabilize $\mathcal{T}\left(  \text{M}\right)  $
and $\Gamma\backslash\mathcal{T}\left(  \text{M}\right)  =\mathfrak{M}%
_{L}\left(  \text{M}\right)  $ is a non-singular quasi-projective variety.
\end{theorem}

\begin{remark}
\label{Vie1}Theorem \ref{Vie} implies that we constructed a family of
non-singular CY manifolds $\pi:\mathcal{X\rightarrow}\mathfrak{M}_{L}\left(
\text{M}\right)  $ over a quasi-projective non-singular variety $\mathfrak{M}%
_{L}($M$)$. Moreover it is easy to see that $\mathcal{X\subset}\mathbb{CP}%
^{N}\times\mathfrak{M}_{L}\left(  \text{M}\right)  .$ \textit{So}
$\mathcal{X}$ \ \textit{is also quasi-projective. From now on we will work
only with this family.}
\end{remark}

\section{Variations of Hodge Structures}

\subsection{Basic Definitions}

\begin{definition}
\label{HS} Let $\Lambda_{g}$ be a free abelian group on which a bilinear form
$\left\langle \text{ , }\right\rangle _{g}$ is defined and $\left\langle
e_{1}\text{,}e_{2}\right\rangle _{n}=(-1)^{n}\left\langle e_{2}\text{,}%
e_{1}\right\rangle _{n}.$ The Hodge Structure of weight $n$ is a decomposition
$\Lambda_{g}\otimes\mathbb{C}=%
{\displaystyle\bigoplus\limits_{p+q=g}}
H^{p,q},$ where $H^{p,q}=\overline{H^{q,p}}$ and the restriction of the
bilinear form $\left\langle \text{ , }\right\rangle $ on $H^{p,q}$ has the
following properties:
\begin{equation}
\left(  -i\right)  ^{p-q}\left(  -1\right)  ^{\frac{g(g-1)}{2}}\left\langle
\omega\text{,}\overline{\omega}\right\rangle _{g}>0, \label{Hod1}%
\end{equation}
For $\omega\in H^{p,q}$ and $\omega_{1}\in H^{p_{1},q_{1}}$ for $p\neq q_{1}$
and $p\neq q_{1}$ we have $\left\langle \omega\text{,}\omega_{1}\right\rangle
_{g}=0.$
\end{definition}

\begin{remark}
\label{Hod2} The Hodge Structure of weight $n$ defines a Hodge filtration on
$\Lambda_{\mathbb{Z}}\otimes\mathbb{C}$:
\begin{equation}
F^{0}\subset F^{1}\subset...\subset F^{g}=\Lambda_{g}\otimes\mathbb{C},
\label{hod1}%
\end{equation}
where $F^{k}=%
{\displaystyle\bigoplus\limits_{p+q=g,q\leq k}}
H^{p,q}.$ The Hodge filtration is isotropic: $\left(  F^{q}\right)  ^{\perp
}=F^{g-q-1}$ and $H^{p,q}=F^{p}\cap\overline{F^{g-p,p}}.$ Once we have an
isotropic filtration, then we have the following Hodge decomposition
$\Lambda_{g}\otimes\mathbb{C}=%
{\displaystyle\bigoplus\limits_{p+q=g}}
H^{p,q}=%
{\displaystyle\bigoplus\limits_{p=0}^{g}}
F^{p}\cap\overline{F^{g-p,p}}.$
\end{remark}

\begin{definition}
\label{VHS}Let $\rho:\pi_{1}(U)\rightarrow Aut\left(  \Lambda_{g}\right)  $ be
a representation of the fundamental group $\pi_{1}$($U$) to $Aut(\Lambda
_{n}).$ Suppose that $G$ preserves the bilinear form $\left\langle \text{ ,
}\right\rangle $. Let $\pi:\widetilde{H}_{\mathbb{C}}=\widetilde{\Lambda_{g}%
}\otimes\mathbb{C\rightarrow}U$ be the flat vector bundle associated with the
above representations of $\pi_{1}(U)$. Let $F_{\tau}^{0}\subset F_{\tau}%
^{1}\subset...\subset F_{\tau}^{g}=\Lambda_{g}\otimes\mathbb{C}$ be an
isotropic filtration which depends holomorphically. We will say that the
representation $\rho$ defines a VHS of weight $n$ over $U$ if the flat
connection $\nabla$ defined by the representation $\rho$ satisfies Griffith's
transversely condition:%
\begin{equation}
\nabla:F_{\tau}^{p}\rightarrow F_{\tau}^{p-1}\otimes\Omega_{U}^{1} \label{GT}%
\end{equation}
and $\Lambda_{g}\otimes\mathbb{C=}%
{\displaystyle\bigoplus\limits_{p=0}^{g}}
F_{\tau}^{p}\cap\overline{F_{\tau}}^{g-q}=%
{\displaystyle\bigoplus\limits_{p+q=g}}
H^{p,q},$ where $H^{p,q}=F^{q}\cap\overline{F}^{g-q}$. We defined the
Variation of Hodge Structure (VHS) of weight $g$ over the complex manifold $U$.
\end{definition}

\subsection{Moduli of VHS}

The moduli spaces $\mathfrak{G}_{g}$ of the VHS of weight $g$ is described as
follows; Let $\mathfrak{G}_{g}($M$)$ be the set of all filtrations $\left(
\ref{hod1}\right)  $ in $\Lambda_{g}\otimes\mathbb{C}$ defined such that
\[
\underset{p+q=g}{\oplus}F^{p}\cap\overline{F}^{g-q}=\underset{p+q=g}{\oplus
}H^{p,q}=\Lambda_{g}\otimes\mathbb{C},
\]
have fixed dimensions $h^{p,q}$ and the filtration is isotropic. Then the
moduli space $\mathfrak{G}_{g}($M$)=\left.  \mathbb{G}\right/  \mathbb{K}%
_{1},$%
\begin{equation}
\mathbb{G}=\left\{
\begin{array}
[c]{c}%
\mathbb{S}p(b_{g})\text{ \ \ \ \ \ \ \ for }g=2m+1\\
\mathbb{SO}(m_{1},m_{2}),\text{ for }g=2m,\text{ }%
\end{array}
\right.  , \label{Hod4}%
\end{equation}
where $m_{1}=%
{\displaystyle\sum\limits_{p=2k\,,\ p\leq m}}
h^{p,q},$ $m_{1}+m_{2}=2m$ and%
\begin{equation}
\mathbb{K}_{1}=\left\{
\begin{array}
[c]{c}%
\left(  \underset{p\leq m}{%
{\displaystyle\prod}
}\mathbb{U}(h^{p,q})\right)  \text{ \ \ \ \ \ \ \ \ \ \ \ \ \ for\ }g=2m+1\\
\left(  \underset{p<m}{%
{\displaystyle\prod}
}\mathbb{U}(h^{p,q})\right)  \times\mathbb{SO}\left(  m\right)  ,\text{ for
}g=2m.\text{ }%
\end{array}
\right.  . \label{mvh}%
\end{equation}
Over $\mathfrak{G}_{g}($M$)=\left.  \mathbb{G}\right/  \mathbb{K}_{1}$ we get
the universal variation of Hodge structure
\begin{equation}
0\subset\mathfrak{F}^{0}\subset...\subset\mathfrak{F}^{g}=H_{g}\times
\mathfrak{G}_{g}=\left(  \Lambda_{g}\otimes\mathbb{C}\right)  \times
\mathfrak{G}_{g} \label{vb}%
\end{equation}
of weight $g$ such that it satisfies $\left(  \ref{GT}\right)  .$

If $\widetilde{H_{g}}=\widetilde{\Lambda_{g}}\otimes\mathbb{C}\rightarrow U$
is a VHS of weight $n$ over $U$, then there exists a map $p:U\rightarrow
\mathfrak{G}_{n}$ such that the VHS $\widetilde{H_{g}}=\widetilde{\Lambda_{g}%
}\otimes\mathbb{C}\rightarrow U$ of weight $g$ over $U$ is the pullback of the
universal VHS $\left(  \ref{vb}\right)  $\ of weight $g$. The map $p$ will be
called the period map.

\begin{remark}
\label{CS}From the definition of $\mathbb{G}$ given by $\left(  \ref{Hod4}%
\right)  $ and the Definition of $\mathbb{K}_{1}$ given by $\left(
\ref{mvh}\right)  $ it is clear that we have $\mathbb{K}_{1}\subseteq
\mathbb{K}$ where $\mathbb{K}$ is a maximal compact subgroup in $\mathbb{G}.$
The compact subvarieties $g\left(  \mathbb{K}/\mathbb{K}_{1}\right)
\subset\mathbb{G}/\mathbb{K}_{1}$ are complex subvarieties for each
$g\in\mathbb{G}$.
\end{remark}

\subsection{$\left.  \mathbb{SU(}g,g)\right/  \mathbb{S(U(}g)\times
\mathbb{U(}g))$ and VHS of Weight $g-$Flat Coordinates}

Let $V:=\mathbb{C}^{2g}$ be equipped with a Hermitian metric
\[
q(z,\overline{z}):=z^{1}\overline{z^{1}}+...+z^{g}\overline{z^{g}}%
-z^{g+1}\overline{z^{g+1}}-...-z^{2g}\overline{z^{2g}}.
\]
The Group that preserve $q(z,\overline{z})$ with determinant $1$ is defined to
be $\mathbb{SU(}g,g)).$ The elements of $\mathbb{SU(}g,g))$ and its maximal
compact subgroup $\mathbb{S(U(}g)\times\mathbb{U(}g))$ are defined as the set
of all $2g\times2g$ matrices respectively:
\begin{equation}
\left(
\begin{array}
[c]{cc}%
A & \left(  \tau_{j}^{i}\right) \\
\left(  \overline{\tau_{j}^{i}}\right)  ^{t} & B
\end{array}
\right)  \text{ and }\left(
\begin{array}
[c]{cc}%
A & 0\\
0 & B
\end{array}
\right)  , \label{2}%
\end{equation}
where $A$ and $B$ are elements of $\mathbb{U(}n)$ such that $\det(A)\det(B)=1$
and the matrices $\left(  \tau_{j}^{i}\right)  $ satisfy $I_{n}-\left(
\tau_{j}^{i}\right)  \times\overline{\left(  \tau_{j}^{i}\right)  }^{t}>0.$ We
define
\[
\mathfrak{h}_{g,g}:=\left.  \mathbb{SU(}g,g)\right/  \mathbb{S(U(}%
g)\times\mathbb{U(}g)).
\]

\begin{theorem}
\label{Sat} Let us consider%
\begin{equation}
\mathcal{D}_{g,g}:=\left\{  \left.  Z=(\tau_{j}^{i})\right\vert 1\leq i,j\leq
g\text{ such that }(I_{g}-Z\times\overline{Z^{t}})>0\right\}  . \label{lc}%
\end{equation}
Then \textbf{A. }$\mathcal{D}_{g,g}=\mathfrak{h}_{g,g}:=\mathbb{SU}%
(g,g)/\mathbb{S}(\mathbb{U}(g)\times\mathbb{U}(g))$ and \textbf{B.
}$\mathfrak{h}_{g,g}\subset\mathbf{Grass}\left(  g,2g\right)  .$
\end{theorem}

\textbf{Proof of A:} Let $H_{2g}$ be a real $2n$ dimensional vector space.
Suppose $E_{2g}:=H_{2g}\otimes\mathbb{C}$ has a Hermitian metric
$q(z,\overline{z})$ of signature $(n,n).$

\begin{lemma}
\label{gr}Let $E_{2g}=E_{g}\oplus\overline{E_{g}}$ such that $q|_{E_{g}}>0$
and $q|_{\overline{E_{g}}}<0.$ Then there is one to one map between
$\mathfrak{h}_{g,g}$ and all subspaces $E_{g}\subset E_{2g}$ of dimension n
such hat $q$%
$\vert$%
$_{E_{g}}>0$.
\end{lemma}

\textbf{Proof: }$\mathbb{SU(}g,g)$ acts transitively of all oriented subspaces
$E_{g}$ of dimension $g$ such that $h|_{E_{g}}>0.$ The stabilizer of a fixed
splitting $E_{g}\oplus\overline{E_{g}}$ such that $q(z,\overline{z})|_{E_{g}%
}>0$ and $q(z,\overline{z})|_{\overline{E_{g}}}<0$ is $\mathbb{S(U(}%
g)\times\mathbb{U(}g)).$ Lemma \ref{gr} is proved. $\blacksquare$ $\ $Part
\textbf{A }of Theorem \ref{Sat} is proved. $\blacksquare$

\textbf{Proof of B: }Theorem \ref{Sat} \textbf{B} follows directly form the
description of the elements of $\mathbb{SU(}n,n))$ and its subgroup
$\mathbb{S(U(}g)\times\mathbb{U(}g))$ given by $\left(  \ref{2}\right)  $ and
since $\mathfrak{h}_{g,g}:=\mathbb{SU}(g,g)/\mathbb{S}(\mathbb{U}%
(g)\times\mathbb{U}(g))$. $\blacksquare$

\begin{theorem}
\label{emb} There exists a representation of\textbf{ }$\mathbb{SU(}g,g)$ to
$\mathbb{S}\mathbf{p}\left(  4g,\mathbb{R}\right)  $ and \textbf{ }%
\[
\mathfrak{h}_{g,g}:=\mathbb{SU}(g,g)/\mathbb{S}(\mathbb{U}(g)\times
\mathbb{U}(g))\subset\mathfrak{h}_{4g}:=\mathbb{S}\mathbf{p}(4g)/\mathbb{U}%
(2g)
\]
is a totally geodesic complex submanifold.
\end{theorem}

\textbf{Proof: }One needs to construct a representation $Sh_{g,g}$ of
$\mathbb{SU(}g,g)$ into a real $4g$ dimensional space $H_{4g}$ such that
$\mathbb{SU(}g,g)$ preserve a skew symmetric non degenerate form $\left(
\text{ },\text{ }\right)  $ $_{4g}$ on $H_{4g}$. The representation $Sh_{g,g}$
defines an embedding
\[
p:\mathfrak{h}_{g,g}=\mathbb{SU(}g,g)/\mathbb{S}\left(  \mathbb{U}\left(
g\right)  \times\mathbb{U}\left(  g\right)  \right)  \subset\mathfrak{h}%
_{4g}=\mathbb{S}\mathbf{p}\mathbb{(}4g)/\mathbb{U}\left(  2g\right)  .
\]
Such a representation is constructed in \cite{Sat}. $\blacksquare$

\begin{definition}
\label{VHSnn1}The symmetric domain $\mathfrak{h}_{2g}=\mathbb{S}%
\mathbf{p}\mathbb{(}4g)/\mathbb{SU}\left(  2g\right)  =\mathfrak{G}_{1}(2g)$
is the moduli space of all VHS of weight one with $h^{1,0}=2g$. Thus the
embedding
\[
p:\mathfrak{h}_{g,g}=\mathbb{SU(}g,g)/\mathbb{S}\left(  \mathbb{U}\left(
g\right)  \times\mathbb{U}\left(  g\right)  \right)  \subset\mathfrak{h}%
_{4g}=\mathbb{S}\mathbf{p}\mathbb{(}4g)/\mathbb{U}\left(  2g\right)  .
\]
induces a VHS of weight one over $\mathfrak{h}_{g,g}.$ Let us denote it by
$\mathcal{V}_{g,g}(1).$ The realization of $\mathcal{V}_{g,g}(1)$ is the
following: Let $H_{2g}\otimes\mathbb{C}=E_{2g}=E_{g}\oplus\overline{E_{g}}$ be
a fixed $2g$ dimensional vector space with the Hermitian bilinear form
$q(z,\overline{z})$ of signature $(n,n)$ such that $E_{g}|_{q(z,\overline{z}%
)}>0$ and $\overline{E_{g}}|_{q(z,\overline{z})}<0$. Let $E_{2g}%
\times\mathfrak{h}_{g,g}\rightarrow\mathfrak{h}_{g,g}$ be the trivial 2n
bundle. The action of $\mathbb{SU(}g,g)$ on $E_{2g}$ and $\mathfrak{h}_{g,g}$
defines a vector bundle $\pi_{g,g}:\mathcal{E}_{g}\rightarrow\mathfrak{h}%
_{g,g}$ where $\pi_{g,g}^{-1}(h)=h\left(  E_{g}\right)  =H_{h(id)}^{1,0}$ and
$\mathcal{E}_{g}\oplus\overline{\mathcal{E}_{g}}=E_{2g}\times\mathfrak{h}%
_{g,g}\rightarrow\mathfrak{h}_{g,g}.$
\end{definition}

\begin{definition}
\label{VHSwn}Let us define VHS $\wedge^{g}\mathcal{V}_{g,g}(1)$ over
$\mathfrak{h}_{g,g}$ as the $g$ exterior power of the VHS of $\mathcal{V}%
_{g,g}(1)$\ of weight one$.$
\end{definition}

\begin{remark}
\label{Bn}The bilinear form $\left\langle \text{ },\text{ }\right\rangle _{n}$
on $\wedge^{g}H_{2g}$ defined by the VHS $\wedge^{g}\mathcal{V}_{g,g}(1)$ over
$\mathfrak{h}_{g,g}$ is the $n$ exterior power of the skew symmetric bilinear
form $\operatorname{Im}q=\left(  \text{ },\text{ }\right)  _{2n}$ on $H_{2g}$
associated with $\mathcal{V}_{g,g}(1)$. Then the symmetric space
\[
\mathfrak{h}_{g,g}:=\left.  \mathbb{SU(}g,g)\right/  \mathbb{S}\left(
\mathbb{U(}g)\times\mathbb{U(}g)\right)
\]
is the moduli space of the VHS $\wedge^{g}\mathcal{V}_{g,g}(1)$ of weight $n.$
\end{remark}

\section{Basic Facts about Jacobians of Hyper-Elliptic Curves}

\subsection{Jacobian of Hyper-Elliptic Curves and the $\theta-$Divisors of
Hyper-Elliptic Curves}

\begin{definition}
\label{Jac}Let $C_{g}$ be any projective algebraic curve of genus $g$. We will
define the Jacobian $J(C_{g})$ of $C_{g}$ as follows:%
\[
J(C_{g}):=\left.  \left\{  \left.  D\right\vert D\text{ is a divisor, }\deg
D=0\right\}  \right/  \backsim,
\]
where $D_{1}\backsim D_{2}$ iff there exists a meromorphic function $f$ on
$C_{g}$ such that $D_{1}-D_{2}=\left(  f\right)  ,$ where $\left(  f\right)  $
is the divisor which consists of all zeroes and all poles with their
multiplicities with a sign minus of the poles of a meromorphic function $f$ on
$C_{g}.$
\end{definition}

Another equivalent definition of $J(C_{g})$ is the following one:

\begin{definition}
\label{JAC}$J(C_{g})$ is the set of all line bundles $\mathcal{O}_{C}\left(
D\right)  ,$ where $D$ is a divisor of degree zero modulo isomorphism.
\end{definition}

It is a well known fact that the line bundle $\mathcal{O}_{C_{g}}\left(
D_{1}\right)  $ is isomorphic to $\mathcal{O}_{C_{g}}\left(  D_{2}\right)  $
if and only if $D_{1}-D_{2}=(f),$ where $f$ is a meromorphic function on
$C_{g}.$

\begin{definition}
\label{Per}Let $C_{g}$ be a curve of genus $g.$ Let $\omega_{1},...,\omega
_{g}$ be a basis of $H^{0}\left(  C,\Omega_{C}^{1}\right)  $ and $\gamma
_{1},...,\gamma_{2g}$ be a basis of $H_{1}\left(  C,\mathbb{Z}\right)  .$ Then
we define the period matrix $\mathcal{P}$
\begin{equation}
\mathcal{P}=\left(
{\displaystyle\int\limits_{\gamma_{i}}}
\omega_{j}\right)  \label{per}%
\end{equation}
with $g$ rows \ and $2g$ columns. It is easy to see that we can choose the
basis $\omega_{1},...,\omega_{g}$ in such way that the matrix $\mathcal{P}$
defined by $\left(  \ref{per}\right)  $\ can be presented as follows:%
\begin{equation}
\left(  E_{g},Z_{g}\right)  \label{perm}%
\end{equation}
where $Z_{g}^{t}=Z_{g},$ i.e. it is a symmetric matrix and $\operatorname{Im}%
Z_{g}$ is positive definite. See \cite{Cl}. Then if $A_{1},...,A_{g}%
;B_{1},...,B_{g}$ is a basis of $H_{1}\left(  C,\mathbb{Z}\right)  $ such
that
\[
\left\langle A_{i},A_{j}\right\rangle =\left\langle B_{i},B_{j}\right\rangle
=0\text{ and }\left\langle A_{i},B_{j}\right\rangle =\delta_{ij}%
\]
then there exists a basis $\left\{  \omega_{i}\right\}  $ of $H^{0}\left(
C_{g},\Omega_{C_{g}}^{1}\right)  $ such that%
\[%
{\displaystyle\int\limits_{A_{j}}}
\omega_{i}=\delta_{ij}\text{ and }%
{\displaystyle\int\limits_{B_{j}}}
\omega_{i}=z_{ij},
\]
where the matrix $\left(  z_{ij}\right)  $ is the symmetric matrix $Z_{g}$
that appears in $\left(  \ref{perm}\right)  .$
\end{definition}

\begin{definition}
There is a natural map $p_{g}:S^{g}\left(  C_{g}\right)  \rightarrow J(C_{g})$
given by%
\begin{equation}
p_{g}\left(  x_{1},...,x_{g}\right)  =\left(
\begin{array}
[c]{ccccc}%
{\displaystyle\int\limits_{\infty}^{x_{1}}}
\omega_{1} & . & . & . &
{\displaystyle\int\limits_{\infty}^{x_{1}}}
\omega_{g}\\
. & . & . & . & .\\
. & . & . & . & .\\
. & . & . & . & .\\%
{\displaystyle\int\limits_{\infty}^{x_{n}}}
\omega_{1} & . & . & . &
{\displaystyle\int\limits_{\infty}^{x_{n}}}
\omega_{g}%
\end{array}
\right)  \in J(C_{g}). \label{f01}%
\end{equation}

\end{definition}

We will need the following Theorems \ref{abel} and \ref{isom}.

\begin{definition}
\label{Pic}Let $C_{g}$ be a curve of genus $g.$ We will define $Pic^{r}\left(
C_{g}\right)  ,$ as the set of all line bundles $\mathcal{O}_{C_{g}}\left(
D\right)  $ such that $\deg D=r,$ modulo isomorphism.
\end{definition}

\begin{theorem}
\label{abel}Let $p_{0}$ be a fixed point on $C_{g}.$ The map
\[
p\rightarrow\mathcal{O}_{C_{g}}(p)\otimes\mathcal{O}_{C_{g}}(-p_{0})
\]
defines an embedding: $C_{g}\subset Pic^{0}(C_{g}).$ Let $p_{1},...,p_{k}%
,q_{1},...,q_{k}$ be points on $C_{g}.$ Then there exists a meromorphic
functions $f$ such that $\left(  f\right)  =%
{\displaystyle\sum\limits_{i=1}^{k}}
p_{i}-%
{\displaystyle\sum\limits_{i=1}^{k}}
q_{i},$ if and only if $%
{\displaystyle\sum\limits_{i=1}^{k}}
\left(  ...%
{\displaystyle\int\limits_{p_{i}}^{q_{i}}}
\omega_{j}...\right)  =0\in\mathbb{C}^{g}.$
\end{theorem}

\subsection{$\theta-$Functions and $\theta-$Divisors in the Jacobians of
Hyper-Elliptic Curves}

\begin{definition}
Let $\delta_{i},$ $\varepsilon_{k}\in\left\{  0,1\right\}  $, $j,$
$k=1,...,g.$ We will define%
\[
\delta=\left[
\begin{array}
[c]{c}%
\delta_{1}\\
.\\
.\\
.\\
\delta_{g}%
\end{array}
\right]  ,\text{ }\varepsilon=\left[
\begin{array}
[c]{c}%
\varepsilon_{1}\\
.\\
.\\
.\\
\varepsilon_{g}%
\end{array}
\right]
\]
and%
\[
\theta%
\genfrac{[}{]}{0pt}{}{\delta}{\varepsilon}%
\left(  z,Z_{g}\right)  =
\]%
\[%
{\displaystyle\sum\limits_{m\in\mathbb{Z}^{g}}}
\exp\left\{  \pi i^{t}\left(  m+\frac{\delta}{2}\right)  Z_{g}\left(
m+\frac{\varepsilon}{2}\right)  +2^{t}\left(  m+\frac{\delta}{2}\right)
\left(  z+\frac{\delta}{2}\right)  \right\}  ,
\]
where $z\in\mathbb{C}^{g}.$ $\theta%
\genfrac{[}{]}{0pt}{}{\delta}{\varepsilon}%
\subset J(C_{g})$ is defined as follows:%
\[
\left\{  z\in J(C_{g})\left\vert \theta%
\genfrac{[}{]}{0pt}{}{\delta}{\varepsilon}%
\left(  z,Z_{g}\right)  =0\right.  \right\}  .
\]
Let $\lambda_{1},...,\lambda_{2g+2}$ be the ramification points of $C_{g}.$ We
will define the $2g+2$ divisors $\theta_{i}$ as the the theta divisor $\theta$
and $2g+1$ divisors $\left(  T_{\lambda_{i}-\infty}\right)  _{\ast}\left(
\theta%
\genfrac{[}{]}{0pt}{}{0}{0}%
\right)  $ obtained by the translation of the $\theta-$divisor $\theta%
\genfrac{[}{]}{0pt}{}{0}{0}%
$\ by the $2g+1$ points of order two $\lambda_{i}-\infty$ on $J(C_{g}).$ It is
clear that
\begin{equation}
\theta_{i}:=\left\{  \left.  \left(  a_{1}+...+a_{g-1}+\lambda_{i}\right)
\right\vert a_{i}\in C_{g}\right\}  . \label{thd}%
\end{equation}

\end{definition}

It is proved by Riemann that $\theta%
\genfrac{[}{]}{0pt}{}{\delta}{\varepsilon}%
$ is a divisor in $J(C_{g}).$ The divisor $\theta%
\genfrac{[}{]}{0pt}{}{0}{0}%
$ is called the $\theta-$divisor in $J(C_{g}).$ See \cite{Cl}, Chapter 4.

We will give another definition of the Theta divisors $\theta_{i}$ based on
Theorem \ref{isom}$.$

\begin{theorem}
\label{isom}$J(C_{g})$ is isomorphic to $Pic^{g-1}(C_{g}).$ See \cite{Cl},
Chapter 4.
\end{theorem}

\begin{definition}
\label{theta}Let $\pi_{g}:C_{g}\rightarrow\mathbb{CP}^{1}$ be a non-singular
hyper-elliptic curve of genus $g\geq3$ given by $y^{2}=F_{2g+1}\left(
z\right)  $ and $P_{g}:S^{g}(C_{g})\rightarrow J(C_{g}).$ Let $\lambda
_{1},...,\lambda_{2g+1}\in\mathbb{C}$ be the zeroes of $F_{2g+1}\left(
z\right)  =0$ and $F_{2g+1}^{\prime}(\lambda_{i})\neq0.$ Then the ramification
points of $\pi_{g}$ are $\lambda_{1},...,\lambda_{2g+1}$ and $\infty.$ We
define the $\theta_{\infty}-$divisor in $J(C_{g})$ as the image of
$S^{g-1}\left(  C_{g}\right)  $ in $Pic^{g-1}(C_{g})$ as follows
\[
\theta_{\infty}:=P_{g}\left(  S^{g-1}\left(  C_{g}\right)  \right)
\otimes\mathcal{O}_{C_{g}}\left(  \infty\right)  \subset J(C_{g}).
\]
We define $\theta_{i}$ be the divisors in $J(C_{g})$ defined as follows
\[
\theta_{i}:=P_{g}\left(  S^{g-1}\left(  C_{g}\right)  \right)  \otimes
\mathcal{O}_{C_{g}}\left(  \lambda_{i}\right)  .
\]

\end{definition}

\begin{remark}
We have%
\begin{equation}
\theta_{i}\cap\theta_{j}:=P_{g}\left(  S^{g-2}\left(  C_{g}\right)  \right)
\otimes\mathcal{O}_{C_{g}}\left(  \lambda_{i}\right)  \otimes\mathcal{O}%
_{C_{g}}\left(  \lambda_{j}\right)  . \label{int1}%
\end{equation}

\end{remark}

\begin{theorem}
\label{bu}Let $C_{g}$ be a hyper-elliptic curve of genus at least two. Let
$S^{g}(C_{g})$ be the symmetric $g^{th}\,\ $power of $C_{g}.$ Then
$S^{g}\left(  C_{g}\right)  $ is obtained from $J\left(  C_{g}\right)  $ by
blowing up the non-singular subvariety $S^{g-2}\left(  C_{g}\right)
\otimes\mathcal{O}_{C_{g}}\left(  2\lambda_{i}\right)  .$
\end{theorem}

\textbf{Proof: }$\pi_{g}:C_{g}\rightarrow\mathbb{CP}^{1}$ is a double cover
ramified over $2g+2$ points $\lambda_{1},...,\lambda_{2g+1}$ and $\infty.$
Thus we get that%
\[
\dim_{\mathbb{C}}H^{0}\left(  C_{g},\mathcal{O}_{C_{g}}\left(  2\lambda
_{i}\right)  \right)  =\dim_{\mathbb{C}}H^{0}\left(  \mathbb{CP}%
^{1},\mathcal{O}_{\mathbb{CP}^{1}}\left(  x_{0}\right)  \right)  =2.
\]
On the other hand we have $2\lambda_{i}\backsim2\lambda_{j}.$ So%
\[%
{\displaystyle\bigcap\limits_{i=1}^{2g+2}}
\theta_{i}=P_{g}\left(  S^{g-2}\left(  C_{g}\right)  \otimes\mathcal{O}%
_{C_{g}}\left(  2\lambda_{i}\right)  \right)  =...=P_{g}\left(  S^{g-2}\left(
C_{g}\right)  \otimes\mathcal{O}_{C_{g}}\left(  2\lambda_{j}\right)  \right)
.
\]
These two facts imply Theorem \ref{bu}. $\blacksquare$

\begin{remark}
\label{int}The definition of the theta divisors and its praimages $P_{g}%
^{-1}\left(  \theta_{i}\right)  $ in $S^{g}\left(  C_{g}\right)  $ implies
that we have%
\begin{equation}
p_{g}^{-1}\left(  \theta_{i}\right)  \cap p_{g}^{-1}\left(  \theta_{j}\right)
=p_{g}^{-1}\left(  \theta_{i}\cap\theta_{j}\right)  =S^{g-2}\left(
C_{g}\right)  \otimes\lambda_{i}\otimes\lambda_{i}. \label{I0}%
\end{equation}

\end{remark}

\begin{definition}
\label{bshc}Let $\widetilde{S^{g}\left(  C_{g}\right)  }$ be obtained by
blowing up first the codimension two submanifold $P_{g}\left(  S^{g-2}\left(
C_{g}\right)  \otimes\mathcal{O}_{C_{g}}\left(  2\lambda_{i}\right)  \right)
=$ $P_{g}\left(  S^{g-2}\left(  C_{g}\right)  \otimes\mathcal{O}_{C_{g}%
}\left(  2\lambda_{j}\right)  \right)  $ in $J(C_{g})$ and then all possible
$p_{g}^{-1}\left(  \theta_{i}\right)  \cap p_{g}^{-1}\left(  \theta
_{j}\right)  $, where%
\[
p_{g}:S^{g}(C_{g})\rightarrow J(C_{g}).
\]

\end{definition}

\begin{remark}
We see that we have the following geometric picture in $\widetilde
{S^{g}\left(  C_{g}\right)  };$ Let $P_{g}:$ $\widetilde{S^{g}\left(
C_{g}\right)  }\rightarrow$ $S^{g}\left(  C_{g}\right)  $ be the holomorphic
blown down map. Then we have $2g+2$ non-singular divisors $P_{g}^{-1}\left(
p_{g}^{-1}\left(  \theta_{i}\right)  \right)  $ intersecting into
$\binom{2g+2}{2}$ non-singular divisors $P_{g}^{-1}\left(  p_{g}^{-1}\left(
\theta_{i}\cap\theta_{j}\right)  \right)  $ and all $2g+2$ divisors are
tangent to the divisor $S^{g-2}\left(  C_{g}\right)  \otimes\mathcal{O}%
_{C_{g}}\left(  2\lambda_{j}\right)  .$ We have exactly the same configuration
as $H_{i},$ $H_{i}\cap H_{j}.$The last two fact implies Lemma \ref{mt}.
$\blacksquare$
\end{remark}

\section{CY\ Manifolds that are Double Covers of $\mathbb{CP}^{g}$ Ramified
over $2g+2$ Hyperplanes in General Position and Symmetric Powers of
Hyper-Elliptic Curves}

\subsection{Basic Construction of CY\ Manifolds that are Double Cover of
$\mathbb{CP}^{g}$ Ramified over $2g+2$ Hyperplanes in General Position}

\begin{definition}
Let $H_{1},...,H_{2g+2}$ be $2g+2$ hyperplanes in $\mathbb{CP}^{g}.$ We will
say that $H_{i}$ are in general position if any different $g$ hyperplanes
intersect transversely in one points and any different $g+1$ do not intersect.
\end{definition}

\begin{theorem}
\label{CYg}Let $H_{1},...,H_{2g+2}$ be $2g+2$ hyperplanes in $\mathbb{CP}^{g}$
in general position. Let us blow up the $\binom{2g+2}{2}$ different
intersection $H_{i}\cap H_{j}$ on $\mathbb{CP}^{g}$ and denote them by
$H_{ij},$ denote the blown up space by $\widetilde{\mathbb{CP}^{g}}.$ Then
there exists a CY manifold that is a double cover M$_{g}\rightarrow
\widetilde{\mathbb{CP}^{g}}$ ramified over the disjoint union of $H_{i}.$
\end{theorem}

\textbf{Proof: }Let us consider the line bundle $\mathcal{O}_{\widetilde
{\mathbb{CP}^{g}}}\left(
{\displaystyle\sum\limits_{i=1}^{2g+2}}
H_{j}\right)  .$ Since as a class of $2g-2$ class of homology is homological
to $12H$ and $\widetilde{\mathbb{CP}^{g}}$ is a simply connected manifold,
then we have%
\[
\mathcal{O}_{\widetilde{\mathbb{CP}^{g}}}\left(
{\displaystyle\sum\limits_{i=1}^{2g+2}}
H_{j}+%
{\displaystyle\sum\limits_{i<j}}
H_{i}\cap H_{j}\right)  \approxeq\mathcal{O}_{\widetilde{\mathbb{CP}^{g}}%
}\left(  12H\right)  .
\]
Now we can define a map of the line bundle
\[
\mathcal{O}_{\widetilde{\mathbb{CP}^{g}}}\left(  6H\right)  =\sqrt
{\mathcal{O}_{\widetilde{\mathbb{CP}^{g}}}\left(
{\displaystyle\sum\limits_{i=1}^{2g+2}}
H_{j}+%
{\displaystyle\sum\limits_{i<j}}
H_{i}\cap H_{j}\right)  }%
\]
to $\mathcal{O}_{\widetilde{\mathbb{CP}^{g}}}\left(
{\displaystyle\sum\limits_{i=1}^{2g+2}}
H_{j}+%
{\displaystyle\sum\limits_{i<j}}
H_{i}\cap H_{j}\right)  $ by taking the map of the fibre $t\rightarrow t^{2}.$
Then the praimage of the section $\sigma$ whose zero set is $%
{\displaystyle\sum\limits_{i=1}^{2g+2}}
H_{j}+%
{\displaystyle\sum\limits_{i<j}}
H_{i}\cap H_{j}$ will define a double cover M$_{n}$ over $\widetilde
{\mathbb{CP}^{g}}$ ramified over all $H_{i}$ and $H_{i}\cap H_{j}.$ Then
direct computation using the adjunction formula shows that M$_{g}$ is a CY
manifold. $\blacksquare$

\subsection{CY\ Manifolds that are Double Covers of $\mathbb{CP}^{g}$ Ramified
over $2g+2$ Hyperplanes in General Position\ and Covered by the $g^{th}$ Power
of A Hyper-Elliptic Curve $C_{g}$ of Genus $g.$}

Let $C_{g}$ be a hyper-elliptic curve of genus $g.$ Then the cover%
\[%
\begin{array}
[c]{ccc}%
\underset{g}{\underbrace{C_{g}\times...\times C_{g}}} & \rightarrow &
S^{g}(C_{g})\\
\downarrow &  & \downarrow\\
\underset{g}{\underbrace{\mathbb{P}^{g}\times...\times\mathbb{P}^{g}}} &
\rightarrow & S^{g}(\mathbb{P}^{g})=\mathbb{P}^{g}%
\end{array}
\]
will be a Galois cover with a Galois group the semi-direct product $\left(
\left.  \mathbb{Z}\right/  2\mathbb{Z}\right)  ^{g}\ltimes S_{g}.$ Let
\begin{equation}
N=\left(  \left.  \mathbb{Z}\right/  2\mathbb{Z}\right)  ^{g-1}\ltimes
S_{g}\subset\left(  \left.  \mathbb{Z}\right/  2\mathbb{Z}\right)  ^{g}\ltimes
S_{g} \label{ge}%
\end{equation}
be subgroup defined as follows:%
\begin{equation}
\left\{  \left(  a_{1,...,}a_{g}\right)  \left\vert a_{i}\in\left.
\mathbb{Z}\right/  2\mathbb{Z}\right.  ,%
{\displaystyle\sum\limits_{i=1}^{g}}
a_{i}=0\right\}  \ltimes S_{g}. \label{ge1}%
\end{equation}

\begin{theorem}
\label{mt}Let $C_{g}$ be a hyper-elliptic curve of genus $g>2.$ Let $N$ be the
subgroup isomorphic $N=\left(  \left.  \mathbb{Z}\right/  2\mathbb{Z}\right)
^{g-1}\ltimes S_{g}$ and embedded into $\left(  \left.  \mathbb{Z}\right/
2\mathbb{Z}\right)  ^{g}\ltimes S_{g}$ by $\left(  \ref{ge1}\right)  .$ Then%
\[
\left.  \underset{g}{\underbrace{C_{g}\times...\times C_{g}}}\right/
N=\text{M}_{g,s}%
\]
is a CY manifold which is a double cover of $\mathbb{CP}^{g}$ ramified over
$2g+2$ hyperplanes in general position. \qquad
\end{theorem}

\textbf{Proof: }Since $N$ is a subgroup of index two in $\left(  \left.
\mathbb{Z}\right/  2\mathbb{Z}\right)  ^{g}\ltimes S_{g}$ and%
\begin{equation}
Q_{g}:\underset{g}{\underbrace{C_{g}\times...\times C_{g}}}\rightarrow\left.
\underset{g}{\underbrace{C_{g}\times...\times C_{g}}}\right/  \left(  \left.
\mathbb{Z}\right/  2\mathbb{Z}\right)  ^{g}\ltimes S_{g}=\mathbb{CP}^{g},
\label{G}%
\end{equation}
be the quotient map. Then there exists a holomorphic double cover%
\begin{equation}
\pi_{g}:\left.  \underset{g}{\underbrace{C_{g}\times...\times C_{g}}}\right/
N=\text{M}_{g,s}\rightarrow\left.  \underset{g}{\underbrace{C_{g}%
\times...\times C_{g}}}\right/  \left(  \left.  \mathbb{Z}\right/
2\mathbb{Z}\right)  ^{g}\ltimes S_{g}=\mathbb{CP}^{g}. \label{cy}%
\end{equation}

\begin{lemma}
\label{mtA}Let us consider
\[
\left.  \left(  \left.  \mathbb{Z}\right/  2\mathbb{Z}\right)  ^{g}\ltimes
S_{g}\right/  \left(  \left.  \mathbb{Z}\right/  2\mathbb{Z}\right)
^{g-1}\ltimes S_{g}=\left.  \mathbb{Z}\right/  2\mathbb{Z}.
\]
Since $\left(  \left.  \mathbb{Z}\right/  2\mathbb{Z}\right)  ^{g}\ltimes
S_{g}$ acts on $\underset{g}{\underbrace{C_{g}\times...\times C_{g}}},$ then
$\left.  \mathbb{Z}\right/  2\mathbb{Z}$ acts on%
\[
\text{M}_{g,s}=\left.  \underset{g}{\underbrace{C_{g}\times...\times C_{g}}%
}\right/  N=\left.  \underset{g}{\underbrace{C_{g}\times...\times C_{g}}%
}\right/  \left(  \left.  \mathbb{Z}\right/  2\mathbb{Z}\right)  ^{g-1}\ltimes
S_{g}.
\]
The ramification divisor $R$ of the action of $\left.  \mathbb{Z}\right/
2\mathbb{Z}$ consists of $2g+2$ hyperplanes in $\mathbb{CP}^{g}$ which are the
images of the $2g+2$ divisors $\lambda_{i}\times\underset{g-1}{\underbrace
{C_{g}\times...\times C_{g}}},$ $i=1,...,2g+2$ under the map%
\[
\underset{g}{\underbrace{C_{g}\times...\times C_{g}}}\rightarrow\left.
\underset{g}{\underbrace{C_{g}\times...\times C_{g}}}\right/  \left(  \left.
\mathbb{Z}\right/  2\mathbb{Z}\right)  ^{g}\ltimes S_{g}=\mathbb{CP}^{g}.
\]

\end{lemma}

\textbf{Proof: }The definition of the action of $N$ and $\left(  \left.
\mathbb{Z}\right/  2\mathbb{Z}\right)  ^{g}\ltimes S_{g}$ on%
\[
\underset{g}{\underbrace{C_{g}\times...\times C_{g}}}%
\]
imply that$\left.  \mathbb{Z}\right/  2\mathbb{Z=}\left.  \left(  \left.
\mathbb{Z}\right/  2\mathbb{Z}\right)  ^{g}\ltimes S_{g}\right/  N$ fixes a
divisor
\[
\left\{  D:=(a_{1},...,a_{g})\right\}  \subset\left.  \underset{g}%
{\underbrace{C_{g}\times...\times C_{g}}}\right/  N
\]
if and only if say $a_{g}=\lambda_{i},$ and $\lambda_{i}$ is a Wierstrass
point on $C_{n}.$ This implies that the ramification divisor of%
\[
\text{M}_{g,s}=\left.  \underset{g}{\underbrace{C_{g}\times...\times C_{g}}%
}\right/  N\rightarrow\left.  \underset{g}{\underbrace{C_{g}\times...\times
C_{g}}}\right/  \left(  \left.  \mathbb{Z}\right/  2\mathbb{Z}\right)
^{g}\ltimes S_{g}=\mathbb{CP}^{g}%
\]
consists of the images of the $2g+2$ divisors $\lambda_{i}\times\underset
{g-1}{\underbrace{C_{g}\times...\times C_{g}}}$ in $\underset{g}%
{\underbrace{C_{g}\times...\times C_{g}}}.$ Since the $2g+2$ points
$\lambda_{i}$ are different in $\mathbb{CP}^{1}$ then the images of
\[
\lambda_{i_{k}}\times\underset{g-1}{\underbrace{C_{g}\times...\times C_{g}}}%
\]
for $k=1,...,g$ and $\lambda_{i}$ different points intersect into the image
$G(\lambda_{i_{1}},...,\lambda_{i_{g}})\in$ $\mathbb{CP}^{g}$ of the point
$(\lambda_{i_{1}},...,\lambda_{i_{g}})$ by the map $\left(  \ref{G}\right)  .$
Any $g+1$ different divisors $\lambda_{i_{k}}\times\underset{g-1}%
{\underbrace{C_{g}\times...\times C_{g}}}$ in $\underset{g}{\underbrace
{C_{g}\times...\times C_{g}}}$ do not intersect. So we get that the
ramification divisor of the double cover $\left(  \ref{cy}\right)  $ is the
union of $2g+2$ hyperplanes that are the images of
\[
\lambda_{i}\times\underset{g-1}{\underbrace{C_{g}\times...\times C_{g}}}%
\]
in $S^{g}\left(  \mathbb{CP}^{1}\right)  =\mathbb{CP}^{g}$ and they are in
general position. So Theorem \ref{CYg} implies Theorem \ref{mt}.
$\blacksquare$

\section{Deformation Theory}

\subsection{Jacobians of Hyper-Elliptic Curves and Locally Symmetric Spaces as
Moduli of CY Manifolds}

\begin{theorem}
\label{KSK}Let $C_{g}$ be a hyper-elliptic curve of genus $g>2$ and
\[
\text{M}_{g,s}=\left.  \underset{g}{\underbrace{C_{g}\times...\times C_{g}}%
}\right/  N
\]
be a CY manifold constructed by Theorem \ref{mt}. Let M$_{g}$ be the minimal
model of M$_{g,s}.$ Then for any $\phi_{1},$ $\phi_{2}\in H^{1}\left(
\text{M}_{g},T_{\text{M}_{g}}^{1,0}\right)  $ we have%
\begin{equation}
\left[  \phi_{1},\phi_{2}\right]  =0. \label{ksk0}%
\end{equation}

\end{theorem}

\textbf{Proof: }Let $J(C_{g})=\left.  \mathbb{C}^{g}\right/  \Lambda.$ If
$(z^{1},...,z^{g})$ are linear coordinates, then
\[
\overline{dz^{1}},...,\overline{dz^{g}}%
\]
will be globally defined anti-holomorphic $(0,1)$ forms on $J(C_{g})=\left.
\mathbb{C}^{g}\right/  \Lambda$ and therefore on $J(C_{g})-\underset{1\leq
i<j\leq2g+2}{\cup}\theta_{ij}.$ By the same arguments we get that%
\[
\frac{\partial}{\partial z^{1}},...,\frac{\partial}{\partial z^{g}}%
\]
are globally defined holomorphic vector field on $J(C_{g})=\left.
\mathbb{C}^{g}\right/  \Lambda$ and therefore on $J(C_{g})-\underset{1\leq
i<j\leq2g+2}{\cup}\theta_{ij}.$

According to Theorem \ref{bu} we have%
\[
J(C_{g})-\underset{1\leq i<j\leq2g+2}{\cup}\theta_{ij}=S^{g}(C_{g}%
)-\underset{1\leq i<j\leq2g+2}{\cup}p_{g}^{-1}\left(  \theta_{j}\cap\theta
_{j}\right)  .
\]
Let%
\[
\Pi:\underset{g}{\underbrace{C_{g}\times...\times C_{g}}}\rightarrow\left.
\underset{g}{\underbrace{C_{g}\times...\times C_{g}}}\right/  S^{g}%
=S^{g}(C_{g}),
\]
be the natural map. Let us consider $\underset{g}{\underbrace{C_{g}%
\times...\times C_{g}}}-R_{g},$ where
\[
R_{g}=\Pi^{-1}\left(  \underset{1\leq i<j\leq2g+2}{\cup}p_{g}^{-1}\left(
\theta_{j}\cap\theta_{j}\right)  \right)  .\
\]
We have a natural etale map%
\begin{equation}
\Phi_{g}:\underset{g}{\underbrace{C_{g}\times...\times C_{g}}}-R_{g}%
\rightarrow\left.  \left(  \underset{g}{\underbrace{C_{g}\times...\times
C_{g}}}-R_{g}\right)  \right/  N, \label{E1}%
\end{equation}
where%
\[
\left.  \left(  \underset{g}{\underbrace{C_{g}\times...\times C_{g}}}%
-R_{g}\right)  \right/  N=\text{M}_{g}-%
{\displaystyle\bigcup\limits_{0<i<j\leq2g+2}}
H_{i}\cap H_{j}%
\]

\begin{lemma}
\label{KSK1}Since the tangent bundle of $J(C_{g})$ is trivial we can choose a
basis of $g$ holomorphic forms $dz^{i}$ and $g$ holomorphic vector fields
$\frac{\partial}{\partial z^{j}}$ on $J(C_{g})$ and thus a basis%
\begin{equation}
\left\{  \overline{dz^{i}}\otimes\frac{\partial}{\partial z^{j}}\right\}  \in
H^{1}\left(  J(C_{g})-\underset{1\leq i<j\leq2g+2}{\cup}\theta_{ij}%
,T_{J(C_{g})-\underset{1\leq i<j\leq2g+2}{\cup}\theta_{ij}}^{1,0}\right)  .
\label{kck4}%
\end{equation}
Let us consider the map:%
\[
P_{g}:\underset{g}{\underbrace{C_{g}\times...\times C_{g}}}-R_{g}\rightarrow
S^{g}(C_{g})-\underset{1\leq i<j\leq2g+2}{\cup}p_{g}^{-1}\left(  \theta
_{j}\cap\theta_{j}\right)  =
\]%
\[
J(C_{g})-\underset{1\leq i<j\leq2g+2}{\cup}\theta_{ij}%
\]
defined by the action of the symmetric group $S_{g}.$ Then there exist a basis
of $g^{2}$ invariant Kodaira-Spencer classes on $\underset{g}{\underbrace
{C_{g}\times...\times C_{g}}}-R_{g},$
\[
\left\{  P_{g}^{\ast}\left(  \overline{dz^{i}}\right)  \otimes P_{g}^{\ast
}\left(  \frac{\partial}{\partial z^{j}}\right)  \right\}  \in H^{1}\left(
\underset{g}{\underbrace{C_{g}\times...\times C_{g}}}-R_{g},T_{\left(
\underset{g}{\underbrace{C_{g}\times...\times C_{g}}}-R_{g}\right)  }%
^{1,0}\right)  ,
\]
invariant under the action of the group $N=\left(  \left.  \mathbb{Z}\right/
2\mathbb{Z}\right)  ^{g-1}\ltimes S_{g}$ defined by $\left(  \ref{ge}\right)
$ and $\left(  \ref{ge1}\right)  .$
\end{lemma}

\textbf{Proof: }Let $\left\{  dz^{i}\right\}  $ be a basis of holomorphic one
forms on $J(C_{g}).$ Let%
\[
P_{g}:\underset{g}{\underbrace{C_{g}\times...\times C_{g}}}\rightarrow
J(C_{g}).
\]
Then we have $P_{g}^{\ast}\left(  \overline{dz^{i}}\right)  =%
{\displaystyle\sum\limits_{k=1}^{g}}
\pi_{k}^{\ast}\left(  \overline{dz^{i}}\right)  ,$ where $\pi_{k}$ is the
projection operator of $\underset{g}{\underbrace{C_{g}\times...\times C_{g}}}$
onto $k$ component $C_{g,k}.$

Let $\left(  P_{g}\right)  _{\ast}\left(  \frac{d}{dz^{j}}\right)  $ be the
dual meromorphic field to to holomorphic one forms $P_{g}^{\ast}\left(
dz^{i}\right)  .$ Since $P_{g}^{\ast}\left(  dz^{i}\right)  $ is an invariant
form under the action of the symmetric group $S_{g}$ and have no zeroes on
$\underset{g}{\underbrace{C_{g}\times...\times C_{g}}}-R_{g},$ then $\left(
P_{g}\right)  _{\ast}\left(  \frac{d}{dz^{j}}\right)  $ are $S_{g}$ invariant
holomorphic vector fields on $\underset{g}{\underbrace{C_{g}\times...\times
C_{g}}}-R_{g}.$ Then
\[
\left\{  P_{g}^{\ast}\left(  \overline{dz^{i}}\right)  \otimes\left(
P_{g}\right)  _{\ast}\left(  \frac{d}{dz^{j}}\right)  \right\}
\]
will be a basis of $g^{2}$ linearly independent Kodaira-Spencer classes on
$\underset{g}{\underbrace{C_{g}\times...\times C_{g}}}-R_{g}$ which are
invariant under the action of the symmetric group $S_{g}$. The group
$N=\left(  \left.  \mathbb{Z}\right/  2\mathbb{Z}\right)  ^{g-1}$ $\ltimes
S_{g}$ by its definitions $\left(  \ref{ge}\right)  $ and $\left(
\ref{ge1}\right)  $ is generated by $g-1$ hyper-elliptic involutions
$\sigma_{g}$ and the symmetric group$.$ The holomorphic $1-$forms $\left\{
dz^{i}\right\}  $ and $\left\{  \frac{\partial}{\partial z^{j}}\right\}  $ on
$S^{g}(C_{g})$ satisfy%
\begin{equation}
\left(  \sigma_{g}\right)  ^{\ast}\left(  \overline{dz^{i}}\right)
=-\overline{dz^{i}}\text{ and }\left(  \sigma_{g}\right)  _{\ast}\left(
\frac{\partial}{\partial z^{j}}\right)  =-\frac{\partial}{\partial z^{j}}.
\label{kck5}%
\end{equation}
Then $\left(  \ref{kck5}\right)  $ implies
\begin{equation}
\left(  \sigma_{g}\right)  ^{\ast}\left(  P_{g}^{\ast}\left(  \overline
{dz^{i}}\right)  \otimes\left(  P_{g}\right)  _{\ast}\left(  \frac{d}{dz^{j}%
}\right)  \right)  =P_{g}^{\ast}\left(  \overline{dz^{i}}\right)
\otimes\left(  P_{g}\right)  _{\ast}\left(  \frac{d}{dz^{j}}\right)  .
\label{IBD}%
\end{equation}
So $\left(  \ref{IBD}\right)  $ implies Lemma \ref{KSK1} . $\blacksquare$

\begin{lemma}
\label{KSK2}Let
\[
\Phi_{g}:\underset{g}{\underbrace{C_{g}\times...\times C_{g}}}-R_{g}%
\rightarrow\left.  \underset{g}{\underbrace{C_{g}\times...\times C_{g}}}%
-R_{g}\right/  N=\text{M}_{g}-\underset{1\leq i<j\leq2g+2}{\cup}H_{i}\cap
H_{j}%
\]
be the quotient map constructed by the action of the group $N$ by Theorem
\ref{mt}. According to Lemma \ref{KSK1} $P_{g}^{\ast}\left(  \overline{dz^{i}%
}\right)  \otimes\left(  P_{g}\right)  _{\ast}\left(  \frac{d}{dz^{j}}\right)
$ are Beltrami differentials on M$_{g}-\underset{1\leq i<j\leq2g+2}{\cup}%
H_{i}\cap H_{j},$ invariant under the action of the group $N$ of index two in
$\left(  \mathbb{Z}/2\mathbb{Z}\right)  ^{g}\ltimes S_{g}.$ Let
\begin{equation}
\left(  \Phi_{g}\right)  _{\ast}\left(  P_{g}^{\ast}\left(  \overline{dz^{i}%
}\right)  \right)  \otimes\left(  \Phi_{g}\right)  _{\ast}\left(  \left(
P_{g}\right)  _{\ast}\left(  \frac{d}{dz^{j}}\right)  \right)  \label{kck6}%
\end{equation}
be the basis of the Kodaira-Spencer classes on M$_{g}-\underset{1\leq
i<j\leq2g+2}{\cup}P_{g}\left(  \theta_{ij}\right)  .$ Then this basis can be
extended to antiholomorphic $(0,1)$ forms$(0,1)$ with coefficients in the
holomorphic tangent bundle on M$_{g}$ and thus as a basis of $H^{1}\left(
\text{M}_{g},T_{\text{M}_{g}}^{1,0}\right)  .$
\end{lemma}

\textbf{Proof: }The $(0,1)$ forms $\overline{dz^{i}}\otimes\frac{\partial
}{\partial z^{j}}$ with coefficients in the holomorphic tangent bundle of
$J(C_{g})-\underset{1\leq i<j\leq2g+2}{\cup}\theta_{ij}$ are invariant under
the holomorphic involution $\sigma_{g}.$ Therefore
\[
\left(  \Phi_{g}\right)  _{\ast}\left(  P_{g}^{\ast}\left(  \overline{dz^{i}%
}\right)  \otimes\left(  P_{g}\right)  _{\ast}\left(  \frac{d}{dz^{j}}\right)
\right)
\]
are well defined Kodaira-Spencer classes in
\[
H^{1}\left(  \text{M}_{g}-\underset{1\leq i<j\leq2g+2}{\cup}H_{i}\cap
H_{j},T_{\text{M}_{g}-\underset{1\leq i<j\leq2g+2}{\cup}H_{i}\cap H_{j}}%
^{1,0}\right)  .
\]
Since the subvarieties $H_{i}\cap H_{j}$ have codimension two, $P_{g}^{\ast
}\left(  \overline{dz^{i}}\right)  \otimes\left(  P_{g}\right)  _{\ast}\left(
\frac{d}{dz^{j}}\right)  $ are invariant under the action of $N.$ Therefore
they are well defined $(0,1)$ forms with coefficients antiholomorphic
functions and holomorphic vector fields. Since they are defined outside
codimension two $\underset{1\leq i<j\leq2g+2}{\cup}P_{g}\left(  \theta
_{ij}\right)  $ in M$_{g},$ then by Hartogs theorem
\[
\left(  \Phi_{g}\right)  _{\ast}\left(  P_{g}^{\ast}\left(  \overline{dz^{i}%
}\right)  \otimes\left(  P_{g}\right)  _{\ast}\left(  \frac{d}{dz^{j}}\right)
\right)
\]
can be extended from M$_{g}-\underset{1\leq i<j\leq2g+2}{\cup}P_{g}\left(
\theta_{ij}\right)  $ to M$_{g}.$ Thus we obtain $g^{2}$ linearly independent
anti-holomorphic forms with coefficients holomorphic vector fields%
\[
\left(  \Phi_{g}\right)  _{\ast}\left(  P_{g}^{\ast}\left(  \overline{dz^{i}%
}\right)  \otimes\left(  P_{g}\right)  _{\ast}\left(  \frac{d}{dz^{j}}\right)
\right)  \in H^{1}\left(  \text{M}_{g},T_{\text{M}_{g}}^{1,0}\right)  .
\]
Since $\dim_{\mathbb{C}}H^{1}\left(  \text{M}_{g},T_{\text{M}_{g}}%
^{1,0}\right)  =g^{2},$ then
\[
\left\{  \left(  \Phi_{g}\right)  _{\ast}\left(  P_{g}^{\ast}\left(
\overline{dz^{i}}\right)  \otimes\left(  P_{g}\right)  _{\ast}\left(  \frac
{d}{dz^{j}}\right)  \right)  \right\}
\]
is a basis of $H^{1}\left(  \text{M}_{g},T_{\text{M}_{g}}^{1,0}\right)  .$
$\blacksquare$

\begin{lemma}
\label{KSK2a}Let $\phi_{k}\in H^{1}\left(  \text{M}_{g},T_{\text{M}_{g}}%
^{1,0}\right)  $ for $k=1$ and $2.$ Then $\left[  \phi_{1},\phi_{2}\right]
=0.$
\end{lemma}

\textbf{Proof: }Since $\left\{  \overline{dz^{i}}\otimes\frac{d}{dz^{j}%
}\right\}  $ is a basis of Kodaira-Spencer classes on the Jacobian
$J_{g}(C_{g}),$ then%
\[
\left\{  P_{g}^{\ast}\left(  \overline{dz^{i}}\right)  \otimes\left(
P_{g}\right)  _{\ast}\left(  \frac{d}{dz^{j}}\right)  \right\}
\]
form a basis of $H^{1}\left(  \text{M}_{g}-\underset{1\leq i<j\leq2g+2}{\cup
}H_{i}\cap H_{j},T_{\text{M}_{g}}^{1,0}\right)  .$ Then any Kodaira-Spencer
classes $\phi_{k},$ $k=1$, $2$ on M$_{g}-\underset{1\leq i<j\leq2g+2}{\cup
}H_{i}\cap H_{j}$ can be expressed as
\begin{equation}
\phi_{k}=%
{\displaystyle\sum\limits_{i,m=1}^{g}}
a_{\overline{i},k}^{m}\left(  \left(  \Phi_{g}\right)  _{\ast}\left(
P_{g}^{\ast}\left(  \overline{dz^{i}}\right)  \otimes\left(  P_{g}\right)
_{\ast}\left(  \frac{d}{dz^{j}}\right)  \right)  \right)  , \label{ksk6}%
\end{equation}
where $a_{\overline{j},k}^{i}$ are constants. Since $H_{i}\cap H_{j}$ have
codimension two, we can extend $\phi_{k}$ on M$_{g}.$ So from $\left(
\ref{ksk6}\right)  ,$ since the coefficients of the $(0,1)$ forms are
antiholomorphic, the vector fields are holomorphic from the definition of
$\left[  \phi_{1},\phi_{2}\right]  $ we get that
\begin{equation}
\left[  \phi_{1},\phi_{2}\right]  =%
{\displaystyle\sum\limits_{\alpha,\beta,m=1}^{g}}
\left(
{\displaystyle\sum\limits_{l=1}^{g}}
\left(  \Phi_{g}\right)  _{\ast}\left(  P_{g}^{\ast}\left(  a_{\overline
{\alpha},1}^{l}\frac{\partial a_{\overline{\beta},2}^{m}}{\partial z^{l}%
}-a_{\overline{\beta},2}^{l}\frac{\partial a_{\overline{a},1}^{m}}{\partial
z^{l}}\right)  \right)  \right)  \overline{dz^{\alpha}}\wedge\overline
{dz^{\beta}}\otimes\frac{\partial}{\partial z^{m}}=0, \label{ksk7}%
\end{equation}
since $a_{\overline{i},k}^{m}$ are constant, then we get that on
M$_{g}-\underset{1\leq i<j\leq2g+2}{\cup}H_{i}\cap H_{j}$
\begin{equation}
\frac{\partial a_{\overline{\beta},2}^{m}}{\partial z^{l}}=\frac{\partial
a_{\overline{a},1}^{m}}{\partial z^{l}}=0. \label{ksk8}%
\end{equation}
Then $\left(  \ref{ksk7}\right)  $ and $\left(  \ref{ksk8}\right)  $ imply
that $\left[  \phi_{1},\phi_{2}\right]  =0$ holds on M$_{g}-\underset{1\leq
i<j\leq2g+2}{\cup}P_{g}\left(  \theta_{ij}\right)  .$ By Hartogs Theorem we
have that $P_{g}^{\ast}\left(  \left[  \phi_{1},\phi_{2}\right]  \right)  =0$
on $J(C_{g})$ and on $\widetilde{S^{g}\left(  C_{g}\right)  }.$ This implies
Lemma \ref{KSK2a}. $\blacksquare$

Lemma \ref{KSK2a} implies Theorem \ref{KSK}. $\blacksquare$

\begin{corollary}
\label{b2}Let M$_{g}$ be a CY manifold which is a double cover of
$\mathbb{CP}^{g}$ for $g\geq3.$ Then the Hodge numbers of M are $h^{g-p,p}%
=\binom{g}{p}^{2}$ and $b_{2}=\binom{2g+2}{2}+1.$
\end{corollary}

\subsection{Moduli space $\mathfrak{M}_{L,g}\left(  \text{M}\right)  $ of
polarized CY manifolds that are double covers of $\mathbb{CP}^{g}$ ramified
over $2g+2$ hyper planes is a Locally Symmetric Spaces}

\begin{theorem}
\label{VPOb}The moduli space $\mathfrak{M}_{L,g}\left(  \text{M}\right)  $ of
polarized CY manifolds that are double covers of $\mathbb{CP}^{g}$ ramified
over $2g+2$ hyper planes is a locally symmetric space $\Gamma_{g}\backslash$
$\mathfrak{h}_{g,g}.$
\end{theorem}

\textbf{Proof: }The proof of Theorem \ref{VPOb} is based on the following
Theorem due to Eli Cartan:

\begin{theorem}
\label{Car}Suppose that X is a quasi projective variety. Suppose that on X
there exists a metric $g$ such that its curvature tensor satisfies $\nabla
R=0.$ Then X is a locally symmetric space.
\end{theorem}

\begin{lemma}
\label{car}The curvature tensor $R$ the Weil-Petersson metric on the moduli
space $\mathfrak{M}_{L,g}\left(  \text{M}\right)  $ of polarized CY manifolds
that are double covers of $\mathbb{CP}^{g}$ ramified over $2g+2$ hyper planes
satisfies $\nabla R=0.$
\end{lemma}

\textbf{Proof: }The proof of Lemma \ref{car} is based on the following Theorem
proved in \cite{to89}:

\begin{theorem}
\label{To} Let $\omega_{\tau}$ be the family of holomorphic forms constructed
by Theorem \ref{tod2} and locally expressed by $\left(  \ref{forms}\right)  .$
Let%
\begin{equation}
\left\Vert \omega_{\tau}\right\Vert _{\mathbf{L}^{2}}^{2}:=(-1)^{\frac
{g(g-1)}{2}}\left(  \frac{\sqrt{-1}}{2}\right)  ^{g}%
{\displaystyle\int\limits_{\text{M}_{\tau}}}
\omega_{\tau}\wedge\overline{\omega_{\tau}}. \label{wp}%
\end{equation}
Then $\log\left(  \left\Vert \omega_{\tau}\right\Vert _{\mathbf{L}^{2}}%
^{2}\right)  $ is a potential of the Weil-Petersson metric, i.e.%
\begin{equation}
dd^{c}\log\left(  \left\Vert \omega_{\tau}\right\Vert _{\mathbf{L}^{2}}%
^{2}\right)  =\operatorname{Im}(W.-P.). \label{wp0}%
\end{equation}

\end{theorem}

We will need the following Proposition, which is reproduced from \cite{to89}:

\begin{proposition}
\label{car0}Let $\left\{  D_{i}^{g}\right\}  $ be a covering of M$_{g},$
$(z_{i}^{1},...,z_{i}^{g})$ be a coordinate system in $D_{i}^{g}$ such that%
\[
\left.  \omega_{\text{M}_{g}}(n,0)\right\vert _{D_{i}^{g}}=dz_{i}^{1}%
\wedge...\wedge dz_{i}^{g}.
\]
According to Lemma \ref{KSK2}
\[
\left\{  \left(  \Phi_{g}\right)  _{\ast}\left(  P_{g}^{\ast}\left(
\overline{dz^{i}}\otimes\frac{\partial}{\partial z^{j}}\right)  \right)
\right\}
\]
is a basis of $H^{1}\left(  \text{M}_{\tau},T_{\text{M}_{\tau}}^{1,0}\right)
.$ Let us consider%
\[
\phi(\tau)=%
{\displaystyle\sum}
\tau_{\beta}^{\alpha}\left(  \left(  \Phi_{g}\right)  _{\ast}\left(
P_{g}^{\ast}\left(  \overline{dz^{\alpha}}\otimes\frac{\partial}{\partial
z^{\beta}}\right)  \right)  \right)  .
\]
Then%
\[
\left(  dz_{i}^{1}+\phi(\tau)\left(  dz_{i}^{1}+\phi(\tau)\right)  \right)
\wedge...\wedge\left(  dz_{i}^{g}+\phi(\tau)\left(  dz_{i}^{g}\right)
\right)  \left.  \omega_{0}\right\vert _{D_{i}^{g}}=
\]%
\begin{equation}
\left.  \omega_{\tau}\right\vert _{D_{i}^{g}}=%
{\displaystyle\sum\limits_{k=1}^{g}}
\left(  -1\right)  ^{\frac{k(k-1)}{2}}%
{\displaystyle\sum\limits_{1\leq i_{l},j_{l}\leq g}}
\tau_{j_{1}}^{i_{1}}\times...\times\tau_{j_{k}}^{i_{k}}\omega_{i_{1}%
j_{1},....,i_{k}j_{k}}\left(  g-k,k\right)  , \label{wp1}%
\end{equation}
where%
\[
\omega_{i_{1}j_{1},....,i_{k}j_{k}}\left(  g-k,k\right)  :=
\]%
\[
\left(  \left(  \Phi_{g}\right)  _{\ast}\left(  P_{g}^{\ast}\left(
\overline{dz^{i_{1}}}\otimes\frac{\partial}{\partial z^{j_{1}}}\right)
\right)  \wedge...\wedge\left(  \Phi_{g}\right)  _{\ast}\left(  P_{g}^{\ast
}\left(  \overline{dz^{i_{k}}}\otimes\frac{\partial}{\partial z^{j_{k}}%
}\right)  \right)  \right)  \lrcorner\omega_{0}.
\]

\end{proposition}

\textbf{Proof: }Corollary \ref{KSK2a} implies that%
\begin{equation}
\phi(\tau):=%
{\displaystyle\sum\limits_{1\leq i,j\leq g}}
\tau_{j}^{i}\left(  \Phi_{g}\right)  _{\ast}\left(  P_{g}^{\ast}\left(
\overline{dz^{i}}\otimes\frac{\partial}{\partial z^{j}}\right)  \right)  .
\label{wp2}%
\end{equation}
Proposition \ref{car0} follows directly from $\left(  \ref{wp2}\right)  ,$
Theorem \ref{tod2} and formula $\left(  \ref{forms}\right)  $ imply $\left(
\ref{wp1}\right)  .$ $\blacksquare$

Repeating the computations of \cite{to89} by using formulas $\left(
\ref{wp}\right)  $ and $\left(  \ref{wp1}\right)  $ we get%
\[
\left\Vert \omega_{\tau}\right\Vert _{\mathbf{L}^{2}}^{2}=\left\Vert
\omega_{0}\right\Vert _{\mathbf{L}^{2}}^{2}+
\]%
\begin{equation}%
{\displaystyle\sum\limits_{1\leq i_{l},j_{l}\leq g}}
\tau_{j_{1}}^{i_{1}}\overline{\tau_{\beta_{1}}^{\alpha_{1}}}\times
...\times\tau_{j_{k}}^{i_{k}}\overline{\tau_{\beta_{k}}^{\alpha_{k}}%
}\left\langle \omega_{i_{1}j_{1},....,i_{k}j_{k}}\left(  g-k,k\right)
,\omega_{\alpha_{1}\beta_{1},....,\alpha_{k}\beta_{k}}\left(  g-k,k\right)
\right\rangle , \label{wp3}%
\end{equation}
where%
\[
\left\langle \omega_{i_{1}j_{1},....,i_{k}j_{k}}\left(  n-k,k\right)
,\omega_{\alpha_{1}\beta_{1},....,\alpha_{k}\beta_{k}}\left(  n-k,k\right)
\right\rangle =
\]%
\begin{equation}
(-1)^{\frac{n(n-1)}{2}}\left(  \frac{\sqrt{-1}}{2}\right)  ^{n}%
{\displaystyle\int\limits_{\text{M}_{0}}}
\omega_{i_{1}j_{1},....,i_{k}j_{k}}\left(  n-k,k\right)  \wedge\overline
{\omega_{\alpha_{1}\beta_{1},....,\alpha_{k}\beta_{k}}\left(  n-k,k\right)  }.
\label{wp4}%
\end{equation}
Formula $\left(  \ref{wp0}\right)  $ implies that the components $g_{i_{k_{1}%
}i_{k_{2}},\overline{j_{k_{1}}}\overline{j_{k_{2}}}}$ of the Weil-Petersson
metric are expressed as follows:%
\[
\frac{\partial^{2}\log\left(  \left\Vert \omega_{\tau}\right\Vert
_{\mathbf{L}^{2}}^{2}\right)  }{\partial\tau_{j_{k1}}^{i_{k_{1}}}%
\overline{\partial\tau_{j_{k2}}^{i_{k_{2}}}}}=
\]%
\begin{equation}
g_{i_{k_{1}}j_{k_{1}},\overline{i_{k_{2}}}\overline{j_{k_{2}}}}=\delta
_{i_{k_{1}}j_{k_{1}},\overline{i_{k_{2}}}\overline{j_{k_{2}}}}-\frac{1}%
{3}R_{i_{k_{1}}j_{k_{1}},\overline{i_{k_{2}}}\overline{j_{k_{2}}}},_{i_{k_{3}%
}j_{k_{3}},\overline{i_{k_{4}}}\overline{j_{k_{4}}}}\tau_{i_{k_{3}}}^{j_{3}%
}\overline{\tau_{i_{k_{4}}}^{j_{4}}}+O(4), \label{wp5}%
\end{equation}
where $R_{i_{k_{1}}j_{k_{1}},\overline{i_{k_{2}}}\overline{j_{k_{2}}}%
},_{i_{k_{3}}j_{k_{3}},\overline{i_{k_{4}}}\overline{j_{k_{4}}}}$ is the
curvature. Formula $\left(  \ref{wp5}\right)  $ implies that in the
coordinates $\tau_{j}^{i},$ the Levi-Cevita covariant derivative of the
Weil-Petersson metric in direction$\frac{\partial}{\partial\tau_{j}^{i}}$ of a
fixed point is given by
\[
\nabla_{\frac{\partial}{\partial\tau_{j}^{i}}}=\frac{\partial}{\partial
\tau_{j}^{i}}.
\]
Since in the expression $\left(  \ref{wp5}\right)  $ we do not have terms of
order three, then%
\begin{equation}
\left.  \nabla_{\frac{\partial}{\partial\tau_{j}^{i}}}R_{i_{k_{1}}i_{k_{2}%
},\overline{j_{k_{1}}}\overline{j_{k_{2}}}},_{i_{k_{3}}i_{k_{3}}%
,\overline{j_{k_{4}}}\overline{j_{k_{4}}}}\right\vert _{\tau=0}=0. \label{wp6}%
\end{equation}
Since the expansion $\left(  \ref{wp5}\right)  $ holds for any point $\tau$ in
$\mathfrak{M}_{L}\left(  \text{M}_{g}\right)  ,$ then we get that $\nabla
R=0.$ Proposition \ref{car0} is proved. $\blacksquare$

Proposition \ref{car0} and Theorem \ref{Car} imply Theorem \ref{VPOb} .
$\blacksquare$

\section{Final Remarks}

\subsection{Some Comments about Miles Reid Conjecture}

M. Reid conjectured that the moduli spaces of polarized CY threefolds is
connected. He suggested that one can connect the CY three-folds by
degenerating them to CY threefold with conic singularity and then deforming
the non-singular CY three-fold obtained by a small resolution of the conic
singularity. It is easy to see that the CY non-singular varieties M$_{g}$ that
are obtained from the double cover of $\mathbb{CP}^{g}$ ramified over $2g+2$
hyper-planes in general position do not contain rational curves with a normal
bundle $\mathcal{O}_{\mathbb{CP}^{1}}\left(  -1\right)  \oplus\mathcal{O}%
_{\mathbb{CP}^{1}}\left(  -1\right)  .$ This implies that the moduli space of
the CY manifold M$_{g}$ is "extremum" in the conjectured connected moduli
space of all CY three-folds. M$_{g}$ is obtained by degeneration of a CY
hypersurfaces in a weighted projective space $\mathbb{CP}^{4}\left(
2:1:...:1\right)  $.

\subsection{The List of Tube Symmetric Domain, which can Parametrize Variation
of Hodge Structure of Weight Three \ Coming from a CY\ Manifold}

In his paper \cite{Gr}, D. Gross classified all symmetric tube domain which
possible can parametrize Variation of Hodge Structure of weight three of
possibly CY threefold. In this paper we realized the Hermitian tube domains
described in \cite{Gr}. It seems that it is a rather difficult problem to
realize the rest of the two tube symmetric domain as Variation of some
geometrically realized CY manifold.

It seems that if the symmetric spaces in the list of Gross can be realized as
Teichm\"{u}ller spaces of CY three-folds, then they should not contain
rational curves with a normal bundles $\mathcal{O}_{\mathbb{CP}^{1}}\left(
-1\right)  \oplus\mathcal{O}_{\mathbb{CP}^{1}}\left(  -1\right)  .$ This
follows from the arithmetic groups that are candidates of the image of the
mapping class groups in the middle cohomologies do not contain elements
$\left(
\begin{array}
[c]{cc}%
1 & 1\\
0 & 1
\end{array}
\right)  ,$ which appears as a monodromy of the conic singularity.

\subsection{Conjectures}

It is easy to see that Weil-Petersson metric on the moduli space of the CY
manifolds M$_{g}$ is a complete metric since it coincides with the Bergman
metric. See also \cite{LSTY}. It is a natural question to ask if the
completeness of the Weil-Petersson metric implies that the moduli space is a
locally symmetric.

Different conjectures about CY manifolds that have locally symmetric spaces as
moduli, can be reviewed \cite{to03}.

\end{document}